\def\temp{1.34}%
\let\tempp=\relax
\expandafter\ifx\csname psboxversion\endcsname\relax
  \message{PSBOX(\temp) loading}%
\else
    \ifdim\temp cm>\psboxversion cm
      \message{PSBOX(\temp) loading}%
    \else
      \message{PSBOX(\psboxversion) is already loaded: I won't load
        PSBOX(\temp)!}%
      \let\temp=\psboxversion
      \let\tempp= 
    \fi
\fi
\tempp
\let\psboxversion=\temp
\catcode`\@=11
%
%
\def\psfortextures{
\def\PSspeci@l##1##2{%
\special{illustration ##1\space scaled ##2}%
}}%
\def\psfordvitops{
\def\PSspeci@l##1##2{%
\special{dvitops: import ##1\space \the\drawingwd \the\drawinght}%
}}%
\def\psfordvips{
\def\PSspeci@l##1##2{%
\d@my=0.1bp \d@mx=\drawingwd \divide\d@mx by\d@my
\includegraphics{##1\space}}}%
\def\psforoztex{
\def\PSspeci@l##1##2{%
\special{##1 \space
      ##2 1000 div dup scale
      \number-\psllx\space \number-\pslly\space translate
}}}%
\def\psfordvitps{
\def\psdimt@n@sp##1{\d@mx=##1\relax\edef\psn@sp{\number\d@mx}}
\def\PSspeci@l##1##2{%
\special{dvitps: Include0 "psfig.psr"}
\psdimt@n@sp{\drawingwd}
\special{dvitps: Literal "\psn@sp\space"}
\psdimt@n@sp{\drawinght}
\special{dvitps: Literal "\psn@sp\space"}
\psdimt@n@sp{\psllx bp}
\special{dvitps: Literal "\psn@sp\space"}
\psdimt@n@sp{\pslly bp}
\special{dvitps: Literal "\psn@sp\space"}
\psdimt@n@sp{\psurx bp}
\special{dvitps: Literal "\psn@sp\space"}
\psdimt@n@sp{\psury bp}
\special{dvitps: Literal "\psn@sp\space startTexFig\space"}
\special{dvitps: Include1 "##1"}
\special{dvitps: Literal "endTexFig\space"}
}}%
\def\psfordvialw{
\def\PSspeci@l##1##2{
\special{language "PostScript",
position = "bottom left",
literal "  \psllx\space \pslly\space translate
  ##2 1000 div dup scale
  -\psllx\space -\pslly\space translate",
include "##1"}
}}%
\def\psforptips{
\def\PSspeci@l##1##2{{
\d@mx=\psurx bp
\advance \d@mx by -\psllx bp
\divide \d@mx by 1000\multiply\d@mx by \xscale
\incm{\d@mx}
\let\tmpx\dimincm
\d@my=\psury bp
\advance \d@my by -\pslly bp
\divide \d@my by 1000\multiply\d@my by \xscale
\incm{\d@my}
\let\tmpy\dimincm
\d@mx=-\psllx bp
\divide \d@mx by 1000\multiply\d@mx by \xscale
\d@my=-\pslly bp
\divide \d@my by 1000\multiply\d@my by \xscale
\at(\d@mx;\d@my){\special{ps:##1 x=\tmpx, y=\tmpy}}
}}}%
\def\psonlyboxes{
\def\PSspeci@l##1##2{%
\at(0cm;0cm){\boxit{\vbox to\drawinght
  {\vss\hbox to\drawingwd{\at(0cm;0cm){\hbox{({\tt##1})}}\hss}}}}
}}%
\def\psloc@lerr#1{%
\let\savedPSspeci@l=\PSspeci@l%
\def\PSspeci@l##1##2{%
\at(0cm;0cm){\boxit{\vbox to\drawinght
  {\vss\hbox to\drawingwd{\at(0cm;0cm){\hbox{({\tt##1}) #1}}\hss}}}}
\let\PSspeci@l=\savedPSspeci@l
}}%
%
%
\newread\pst@mpin
\newdimen\drawinght\newdimen\drawingwd
\newdimen\psxoffset\newdimen\psyoffset
\newbox\drawingBox
\newcount\xscale \newcount\yscale \newdimen\pscm\pscm=1cm
\newdimen\d@mx \newdimen\d@my
\newdimen\pswdincr \newdimen\pshtincr
\let\ps@nnotation=\relax
{\catcode`\|=0 |catcode`|\=12 |catcode`|
|catcode`#=12 |catcode`*=14
|xdef|backslashother{\}*
|xdef|percentother{
|xdef|tildeother{~}*
|xdef|sharpother{#}*
}%
\def\R@moveMeaningHeader#1:->{}%
\def\uncatcode#1{%
\edef#1{\expandafter\R@moveMeaningHeader\meaning#1}}%
\def\execute#1{#1}
\def\psm@keother#1{\catcode`#112\relax}
\def\executeinspecs#1{%
\execute{\begingroup\let\do\psm@keother\dospecials\catcode`\^^M=9#1\endgroup}}%
\def\@mpty{}%
\def\matchexpin#1#2{
  \fi%
  \edef\tmpb{{#2}}%
  \expandafter\makem@tchtmp\tmpb%
  \edef\tmpa{#1}\edef\tmpb{#2}%
  \expandafter\expandafter\expandafter\m@tchtmp\expandafter\tmpa\tmpb\endm@tch%
  \if\match%
}%
\def\matchin#1#2{%
  \fi%
  \makem@tchtmp{#2}%
  \m@tchtmp#1#2\endm@tch%
  \if\match%
}%
\def\makem@tchtmp#1{\def\m@tchtmp##1#1##2\endm@tch{%
  \def\tmpa{##1}\def\tmpb{##2}\let\m@tchtmp=\relax%
  \ifx\tmpb\@mpty\def\match{YN}%
  \else\def\match{YY}\fi%
}}%
\def\incm#1{{\psxoffset=1cm\d@my=#1
 \d@mx=\d@my
  \divide\d@mx by \psxoffset
  \xdef\dimincm{\number\d@mx.}
  \advance\d@my by -\number\d@mx cm
  \multiply\d@my by 100
 \d@mx=\d@my
  \divide\d@mx by \psxoffset
  \edef\dimincm{\dimincm\number\d@mx}
  \advance\d@my by -\number\d@mx cm
  \multiply\d@my by 100
 \d@mx=\d@my
  \divide\d@mx by \psxoffset
  \xdef\dimincm{\dimincm\number\d@mx}
}}%
%
\newif\ifNotB@undingBox
\newhelp\PShelp{Proceed: you'll have a 5cm square blank box instead of
your graphics (Jean Orloff).}%
\def\s@tsize#1 #2 #3 #4\@ndsize{
  \def\psllx{#1}\def\pslly{#2}%
  \def\psurx{#3}\def\psury{#4}
  \ifx\psurx\@mpty\NotB@undingBoxtrue
  \else
    \drawinght=#4bp\advance\drawinght by-#2bp
    \drawingwd=#3bp\advance\drawingwd by-#1bp
  \fi
  }%
\def\sc@nBBline#1:#2\@ndBBline{\edef\p@rameter{#1}\edef\v@lue{#2}}%
\def\g@bblefirstblank#1#2:{\ifx#1 \else#1\fi#2}%
{\catcode`\%=12
\xdef\B@undingBox{
\def\ReadPSize#1{
 \readfilename#1\relax
 \let\PSfilename=\lastreadfilename
 \openin\pst@mpin=#1\relax
 \ifeof\pst@mpin \errhelp=\PShelp
   \errmessage{I haven't found your postscript file (\PSfilename)}%
   \psloc@lerr{was not found}%
   \s@tsize 0 0 142 142\@ndsize
   \closein\pst@mpin
 \else
   \if\matchexpin{\GlobalInputList}{, \lastreadfilename}%
   \else\xdef\GlobalInputList{\GlobalInputList, \lastreadfilename}%
     \immediate\write\psbj@inaux{\lastreadfilename,}%
   \fi%
   \loop
     \executeinspecs{\catcode`\ =10\global\read\pst@mpin to\n@xtline}%
     \ifeof\pst@mpin
       \errhelp=\PShelp
       \errmessage{(\PSfilename) is not an Encapsulated PostScript File:
           I could not find any \B@undingBox: line.}%
       \edef\v@lue{0 0 142 142:}%
       \psloc@lerr{is not an EPSFile}%
       \NotB@undingBoxfalse
     \else
       \expandafter\sc@nBBline\n@xtline:\@ndBBline
       \ifx\p@rameter\B@undingBox\NotB@undingBoxfalse
         \edef\t@mp{%
           \expandafter\g@bblefirstblank\v@lue\space\space\space}%
         \expandafter\s@tsize\t@mp\@ndsize
       \else\NotB@undingBoxtrue
       \fi
     \fi
   \ifNotB@undingBox\repeat
   \closein\pst@mpin
 \fi
\message{#1}%
}%
%
%
\def\psboxto(#1;#2)#3{\vbox{%
   \ReadPSize{#3}%
   \advance\pswdincr by \drawingwd
   \advance\pshtincr by \drawinght
   \divide\pswdincr by 1000
   \divide\pshtincr by 1000
   \d@mx=#1
   \ifdim\d@mx=0pt\xscale=1000
         \else \xscale=\d@mx \divide \xscale by \pswdincr\fi
   \d@my=#2
   \ifdim\d@my=0pt\yscale=1000
         \else \yscale=\d@my \divide \yscale by \pshtincr\fi
   \ifnum\yscale=1000
         \else\ifnum\xscale=1000\xscale=\yscale
                    \else\ifnum\yscale<\xscale\xscale=\yscale\fi
              \fi
   \fi
   \divide\drawingwd by1000 \multiply\drawingwd by\xscale
   \divide\drawinght by1000 \multiply\drawinght by\xscale
   \divide\psxoffset by1000 \multiply\psxoffset by\xscale
   \divide\psyoffset by1000 \multiply\psyoffset by\xscale
   \global\divide\pscm by 1000
   \global\multiply\pscm by\xscale
   \multiply\pswdincr by\xscale \multiply\pshtincr by\xscale
   \ifdim\d@mx=0pt\d@mx=\pswdincr\fi
   \ifdim\d@my=0pt\d@my=\pshtincr\fi
   \message{scaled \the\xscale}%
 \hbox to\d@mx{\hss\vbox to\d@my{\vss
   \global\setbox\drawingBox=\hbox to 0pt{\kern\psxoffset\vbox to 0pt{%
      \kern-\psyoffset
      \PSspeci@l{\PSfilename}{\the\xscale}%
      \vss}\hss\ps@nnotation}%
   \global\wd\drawingBox=\the\pswdincr
   \global\ht\drawingBox=\the\pshtincr
   \global\drawingwd=\pswdincr
   \global\drawinght=\pshtincr
   \baselineskip=0pt
   \copy\drawingBox
 \vss}\hss}%
  \global\psxoffset=0pt
  \global\psyoffset=0pt
  \global\pswdincr=0pt
  \global\pshtincr=0pt 
  \global\pscm=1cm 
}}%
%
%
\def\psboxscaled#1#2{\vbox{%
  \ReadPSize{#2}%
  \xscale=#1
  \message{scaled \the\xscale}%
  \divide\pswdincr by 1000 \multiply\pswdincr by \xscale
  \divide\pshtincr by 1000 \multiply\pshtincr by \xscale
  \divide\psxoffset by1000 \multiply\psxoffset by\xscale
  \divide\psyoffset by1000 \multiply\psyoffset by\xscale
  \divide\drawingwd by1000 \multiply\drawingwd by\xscale
  \divide\drawinght by1000 \multiply\drawinght by\xscale
  \global\divide\pscm by 1000
  \global\multiply\pscm by\xscale
  \global\setbox\drawingBox=\hbox to 0pt{\kern\psxoffset\vbox to 0pt{%
     \kern-\psyoffset
     \PSspeci@l{\PSfilename}{\the\xscale}%
     \vss}\hss\ps@nnotation}%
  \advance\pswdincr by \drawingwd
  \advance\pshtincr by \drawinght
  \global\wd\drawingBox=\the\pswdincr
  \global\ht\drawingBox=\the\pshtincr
  \global\drawingwd=\pswdincr
  \global\drawinght=\pshtincr
  \baselineskip=0pt
  \copy\drawingBox
  \global\psxoffset=0pt
  \global\psyoffset=0pt
  \global\pswdincr=0pt
  \global\pshtincr=0pt 
  \global\pscm=1cm
}}%
%
\def\psbox#1{\psboxscaled{1000}{#1}}%
\newif\ifn@teof\n@teoftrue
\newif\ifc@ntrolline
\newif\ifmatch
\newread\j@insplitin
\newwrite\j@insplitout
\newwrite\psbj@inaux
\immediate\openout\psbj@inaux=psbjoin.aux
\immediate\write\psbj@inaux{\string\joinfiles}%
\immediate\write\psbj@inaux{\jobname,}%
%
%
\def\toother#1{\ifcat\relax#1\else\expandafter%
  \toother@ux\meaning#1\endtoother@ux\fi}%
\def\toother@ux#1 #2#3\endtoother@ux{\def\tmp{#3}%
  \ifx\tmp\@mpty\def\tmp{#2}\let\next=\relax%
  \else\def\next{\toother@ux#2#3\endtoother@ux}\fi%
\next}%
%
%
\let\readfilenamehook=\relax
\def\re@d{\expandafter\re@daux}
\def\re@daux{\futurelet\nextchar\stopre@dtest}%
\def\re@dnext{\xdef\lastreadfilename{\lastreadfilename\nextchar}%
  \afterassignment\re@d\let\nextchar}%
\def\stopre@d{\egroup\readfilenamehook}%
\def\stopre@dtest{%
  \ifcat\nextchar\relax\let\nextread\stopre@d
  \else
    \ifcat\nextchar\space\def\nextread{%
      \afterassignment\stopre@d\chardef\nextchar=`}%
    \else\let\nextread=\re@dnext
      \toother\nextchar
      \edef\nextchar{\tmp}%
    \fi
  \fi\nextread}%
\def\readfilename{\bgroup%
  \let\\=\backslashother \let\%=\percentother \let\~=\tildeother
  \let\#=\sharpother \xdef\lastreadfilename{}%
  \re@d}%
%
%
\xdef\GlobalInputList{\jobname}%
\def\psnewinput{%
  \def\readfilenamehook{
    \if\matchexpin{\GlobalInputList}{, \lastreadfilename}%
    \else\xdef\GlobalInputList{\GlobalInputList, \lastreadfilename}%
      \immediate\write\psbj@inaux{\lastreadfilename,}%
    \fi%
    \ps@ldinput\lastreadfilename\relax%
    \let\readfilenamehook=\relax%
  }\readfilename%
}%
\expandafter\ifx\csname @@input\endcsname\relax    
  \immediate\let\ps@ldinput=\input\def\input{\psnewinput}%
\else
  \immediate\let\ps@ldinput=\@@input
  \def\@@input{\psnewinput}%
\fi%
\def\nowarnopenout{%
 \def\warnopenout##1##2{%
   \readfilename##2\relax
   \message{\lastreadfilename}%
   \immediate\openout##1=\lastreadfilename\relax}}%
\def\warnopenout#1#2{%
 \readfilename#2\relax
 \def\t@mp{TrashMe,psbjoin.aux,psbjoint.tex,}\uncatcode\t@mp
 \if\matchexpin{\t@mp}{\lastreadfilename,}%
 \else
   \immediate\openin\pst@mpin=\lastreadfilename\relax
   \ifeof\pst@mpin
     \else
     \errhelp{If the content of this file is so precious to you, abort (ie
press x or e) and rename it before retrying.}%
     \errmessage{I'm just about to replace your file named \lastreadfilename}%
   \fi
   \immediate\closein\pst@mpin
 \fi
 \message{\lastreadfilename}%
 \immediate\openout#1=\lastreadfilename\relax}%
{\catcode`\%=12\catcode`\*=14
\gdef\splitfile#1{*
 \readfilename#1\relax
 \immediate\openin\j@insplitin=\lastreadfilename\relax
 \ifeof\j@insplitin
   \message{! I couldn't find and split \lastreadfilename!}*
 \else
   \immediate\openout\j@insplitout=TrashMe
   \message{< Splitting \lastreadfilename\space into}*
   \loop
     \ifeof\j@insplitin
       \immediate\closein\j@insplitin\n@teoffalse
     \else
       \n@teoftrue
       \executeinspecs{\global\read\j@insplitin to\spl@tinline\expandafter
         \ch@ckbeginnewfile\spl@tinline
       \ifc@ntrolline
       \else
         \toks0=\expandafter{\spl@tinline}*
         \immediate\write\j@insplitout{\the\toks0}*
       \fi
     \fi
   \ifn@teof\repeat
   \immediate\closeout\j@insplitout
 \fi\message{>}*
}*
\gdef\ch@ckbeginnewfile#1
 \def\t@mp{#1}*
 \ifx\@mpty\t@mp
   \def\t@mp{#3}*
   \ifx\@mpty\t@mp
     \global\c@ntrollinefalse
   \else
     \immediate\closeout\j@insplitout
     \warnopenout\j@insplitout{#2}*
     \global\c@ntrollinetrue
   \fi
 \else
   \global\c@ntrollinefalse
 \fi}*
\gdef\joinfiles#1\into#2{*
 \message{< Joining following files into}*
 \warnopenout\j@insplitout{#2}*
 \message{:}*
 {*
 \edef\w@##1{\immediate\write\j@insplitout{##1}}*
\w@{
\w@{
\w@{
\w@{
\w@{
\w@{
\w@{
\w@{
\w@{
\w@{
\w@{\string\input\space psbox.tex}*
\w@{\string\splitfile{\string\jobname}}*
\w@{\string\let\string\autojoin=\string\relax}*
}*
 \expandafter\tre@tfilelist#1, \endtre@t
 \immediate\closeout\j@insplitout
 \message{>}*
}*
\gdef\tre@tfilelist#1, #2\endtre@t{*
 \readfilename#1\relax
 \ifx\@mpty\lastreadfilename
 \else
   \immediate\openin\j@insplitin=\lastreadfilename\relax
   \ifeof\j@insplitin
     \errmessage{I couldn't find file \lastreadfilename}*
   \else
     \message{\lastreadfilename}*
     \immediate\write\j@insplitout{
     \executeinspecs{\global\read\j@insplitin to\oldj@ininline}*
     \loop
       \ifeof\j@insplitin\immediate\closein\j@insplitin\n@teoffalse
       \else\n@teoftrue
         \executeinspecs{\global\read\j@insplitin to\j@ininline}*
         \toks0=\expandafter{\oldj@ininline}*
         \let\oldj@ininline=\j@ininline
         \immediate\write\j@insplitout{\the\toks0}*
       \fi
     \ifn@teof
     \repeat
   \immediate\closein\j@insplitin
   \fi
   \tre@tfilelist#2, \endtre@t
 \fi}*
}%
\def\autojoin{%
 \immediate\write\psbj@inaux{\string\into{psbjoint.tex}}%
 \immediate\closeout\psbj@inaux
 \expandafter\joinfiles\GlobalInputList\into{psbjoint.tex}%
}%
%
%
%
\def\centinsert#1{\midinsert\line{\hss#1\hss}\endinsert}%
\def\psannotate#1#2{\vbox{%
  \def\ps@nnotation{#2\global\let\ps@nnotation=\relax}#1}}%
\def\pscaption#1#2{\vbox{%
   \setbox\drawingBox=#1
   \copy\drawingBox
   \vskip\baselineskip
   \vbox{\hsize=\wd\drawingBox\setbox0=\hbox{#2}%
     \ifdim\wd0>\hsize
       \noindent\unhbox0\tolerance=5000
    \else\centerline{\box0}%
    \fi
}}}%
%
\def\at(#1;#2)#3{\setbox0=\hbox{#3}\ht0=0pt\dp0=0pt
  \rlap{\kern#1\vbox to0pt{\kern-#2\box0\vss}}}%
%
\newdimen\gridht \newdimen\gridwd
\def\gridfill(#1;#2){%
  \setbox0=\hbox to 1\pscm
  {\vrule height1\pscm width.4pt\leaders\hrule\hfill}%
  \gridht=#1
  \divide\gridht by \ht0
  \multiply\gridht by \ht0
  \gridwd=#2
  \divide\gridwd by \wd0
  \multiply\gridwd by \wd0
  \advance \gridwd by \wd0
  \vbox to \gridht{\leaders\hbox to\gridwd{\leaders\box0\hfill}\vfill}}%
%
\def\fillinggrid{\at(0cm;0cm){\vbox{%
  \gridfill(\drawinght;\drawingwd)}}}%
%
%
\def\textleftof#1:{%
  \setbox1=#1
  \setbox0=\vbox\bgroup
    \advance\hsize by -\wd1 \advance\hsize by -2em}%
\def\textrightof#1:{%
  \setbox0=#1
  \setbox1=\vbox\bgroup
    \advance\hsize by -\wd0 \advance\hsize by -2em}%
\def\endtext{%
  \egroup
  \hbox to \hsize{\valign{\vfil##\vfil\cr%
\box0\cr%
\noalign{\hss}\box1\cr}}}%
%
\def\frameit#1#2#3{\hbox{\vrule width#1\vbox{%
  \hrule height#1\vskip#2\hbox{\hskip#2\vbox{#3}\hskip#2}%
        \vskip#2\hrule height#1}\vrule width#1}}%
\def\boxit#1{\frameit{0.4pt}{0pt}{#1}}%
\catcode`\@=12 
%
 \psfordvips   

\magnification=\magstep1
\font\Bbb=msbm10
\font\Bbbs=msbm8
\textfont12=\Bbb
\scriptfont12=\Bbbs
\font\rmsmall= cmr8
\font\bfone= cmbx10 scaled \magstep1
\let\mcd=\mathchardef
\mcd\Be="7C42
\mcd\Ee="7C45
\mcd\Je="7C4A
\mcd\Pe="7C50
\mcd\Re="7C52
\mcd\Ze="7C5A
\def\al{\alpha}
\def\ep{\epsilon}
\def\la{\lambda}
\def\bra{\langle}
\def\ket{\rangle}
\def\tih{\tilde h}
\def\bah{\bar h}
\def\S{\cal S}
\centerline{\bfone A necklace of Wulff shapes}
\bigskip
\centerline{Jo\"el De Coninck\footnote{$^{(1)}$}
{\rmsmall Centre de Recherche en Mod\'elisation Mol\'eculaire, Universit\'e de 
Mons-Hainaut, 20 Place du Parc, 7000 Mons, Belgium. 
Email: Joel.De.Coninck@crmm.umh.ac.be}, 
Fran\c cois Dunlop\footnote{$^{(2)}$}
{\rmsmall Laboratoire de Physique Th\'eorique et Mod\'elisation (CNRS - UMR 
8089), Universit\'e de Cergy-Pontoise, 95302 Cergy-Pontoise, France. 
Email: dunlop@ptm.u-cergy.fr, huillet@ptm.u-cergy.fr}, 
Thierry Huillet$^{(2)}$}
\bigskip\noindent{\bf Abstract:} {\rmsmall
In a probabilistic model of a film over a disordered substrate, Monte-Carlo
simulations show that the film hangs from peaks of the substrate. The film
profile is well approximated by a necklace of Wulff shapes. Such a necklace
can be obtained as the infimum of a collection of Wulff shapes resting on the 
substrate. When the random substrate is given by iid heights with exponential 
distribution, we prove estimates on the probability density of the resulting
peaks, at small density.}
 
\medskip\noindent {\rmsmall KEYWORDS: Interfaces, random substrate, Wulff 
shape, SOS model}
 
\bigskip\noindent
{\bf 1. Introduction}
\medskip\noindent
A problem in the science of coating is to characterize the surface of a coated 
material as function of the substrate surface and properties of the coating 
material (polymer, resin, metal\dots). 
The topography of a substrate has an important influence on the properties of 
the considered material in terms of lubrication, optical properties, wetting...
Moreover, it is often desirable to coat this substrate with a thin film to 
protect the material. Of course, the topography of the thin film surface and of
the substrate do not have to be the same. This will be a function of the film 
thickness.

A statistical mechanical model of a film requires at least two parameters, one
associated with a surface tension, the other with a pressure difference or
chemical potential controlling the film thickness. The version of the 
Solid-On-Solid model introduced by Abraham and Smith [AS1, AS2] is the simplest
such model. In Section 2, we use it first to model a disordered substrate and
then, with a different set of parameters, to model a film on top of the 
substrate. Numerical simulations show that typical configurations look like a
necklace of Wulff shapes suspended from the peaks of the substrate. This 
motivates Section 3, where only the substrate is random, and the film is 
defined as the infimum of a collection of Wulff shapes over the substrate.
In Section 4 we prove estimates on the density of relevant substrate peaks,
when the pressure difference goes to zero. In section 5 we give a Gibbs 
formulation for the probability density of substrate peak localizations and 
heights.
\bigskip\noindent
{\bf 2. Solid-On-Solid film over quenched Solid-On-Solid substrate}
\medskip\noindent
The substrate is a one-dimensional Solid-On-Solid model with Hamitonian
$$
H_1=J_1\sum_{|i-j|=1}|h^1_i-h^1_j|+K_1\sum h^1_i \eqno{(2.1)}
$$
where $h^1_i\in\Re_+$ is the height of the substrate surface at point 
$i\in\Ze$. It can also represent the height of a first coating, fixed before
the later deposit of a film. 

The substrate being generated according to the Gibbs measure with Hamiltonian 
(2.1) at temperature $kT=1$, and then quenched, a film is deposited and 
thermalized according to the Gibbs measure with Hamiltonian
$$
H_2=J_2\sum_{|i-j|=1}|h^2_i-h^2_j|+K_2\sum h^2_i \eqno{(2.2)}
$$
at temperature $kT=1$. The film height $h^2_i\in\Re_+$ at point $i$ obeys
$h^2_i\ge h^1_i$. This is a grand canonical ensemble where the film volume is 
controlled by $K_2$ while being allowed to fluctuate. In the thermodynamic 
limit, the properties will be the same as if obtained from a fixed volume 
ensemble where $\sum (h^2_i-h^1_i)$ is fixed.

The resulting model of a film on a quenched substrate is studied by Monte-Carlo
simulation with periodic boundary conditions and a heat bath algorithm. Fig. 1 
shows one substrate $h^1$ and, on top, thermal averaged films $\bar h^2$ at 
various values of the pressure $K_2$. 
\vskip-6cm
$$
\psboxto(13cm;0cm){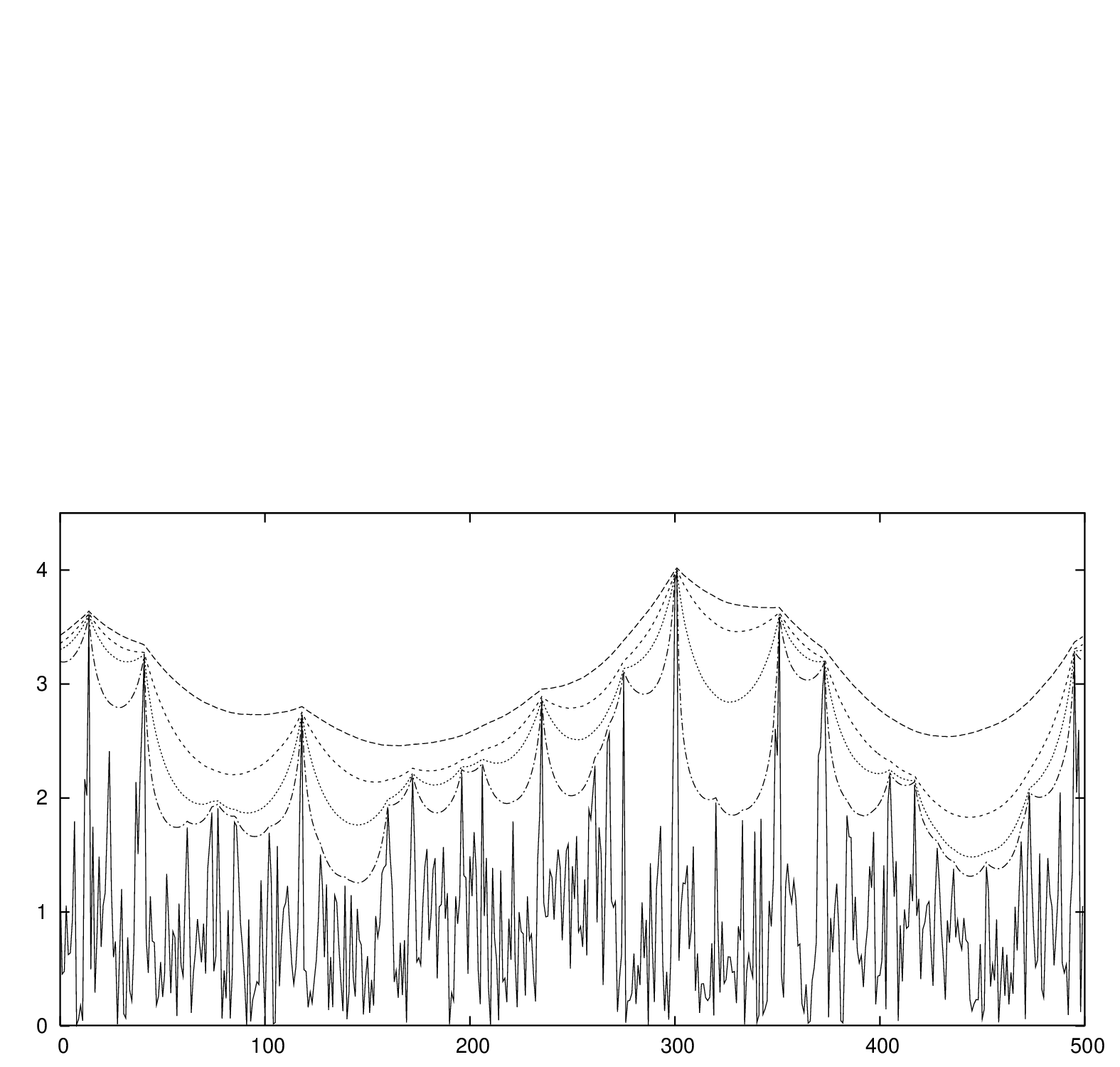}
$$
\vskip-0.2cm\centerline
{\rmsmall Fig. 1: Substrate $J_1=1,K_1=0.5$; film $J_2=30$ and from top to 
bottom $K_2=0.2,\,0.5,\,1,\,2$.}
\bigskip
For small enough $K_2/J_2$, the film appears to hang as from a set of 
telegraph poles of random heights. Between two successive poles, the film 
profile can be checked to be very near a Wulff shape associated with (2.2) at 
the corresponding value of $K_2$. Indeed Fig. 2 shows a portion of substrate, 
a thermal averaged film with $J_2=5$ and $K_2=0.125$, and a necklace of Wulff 
shapes: Each piece is a translate of one and same Wulff shape, with the same 
$J_2$ and $K_2$ as the film. 
The parametric equations of the Wulff shape are [BN,DD]: 
$$
\eqalign{
x(\tan\theta)&=-{1\over K_2}
{d\over d\bigl(\tan\theta\bigr)}\tilde\sigma\bigl(\tan\theta\bigr)\cr
z(\tan\theta)&=-{1\over K_2}\Bigl(\tilde\sigma\bigl(\tan\theta\bigr)-\tan\theta
{d\over d\bigl(\tan\theta\bigr)}\tilde\sigma\bigl(\tan\theta\bigr)\Bigr)
} \eqno{(2.3)}
$$
where $\tilde\sigma(\tan\theta)$ is the projected surface tension, which for
the Solid-On-Solid model takes the form [DD]
$$
\tilde\sigma(\tan\theta)=f(\tan\theta)
-\log\Bigl({f(\tan\theta)+2\over J_2}\Bigr) \eqno{(2.4)}
$$
with
$$
f(\tan\theta)=\bigl(1+(J_2\tan\theta)^2\bigr)^{1/2}-1 \eqno{(2.5)}
$$
\vskip-6cm
$$
\psboxto(0cm;12cm){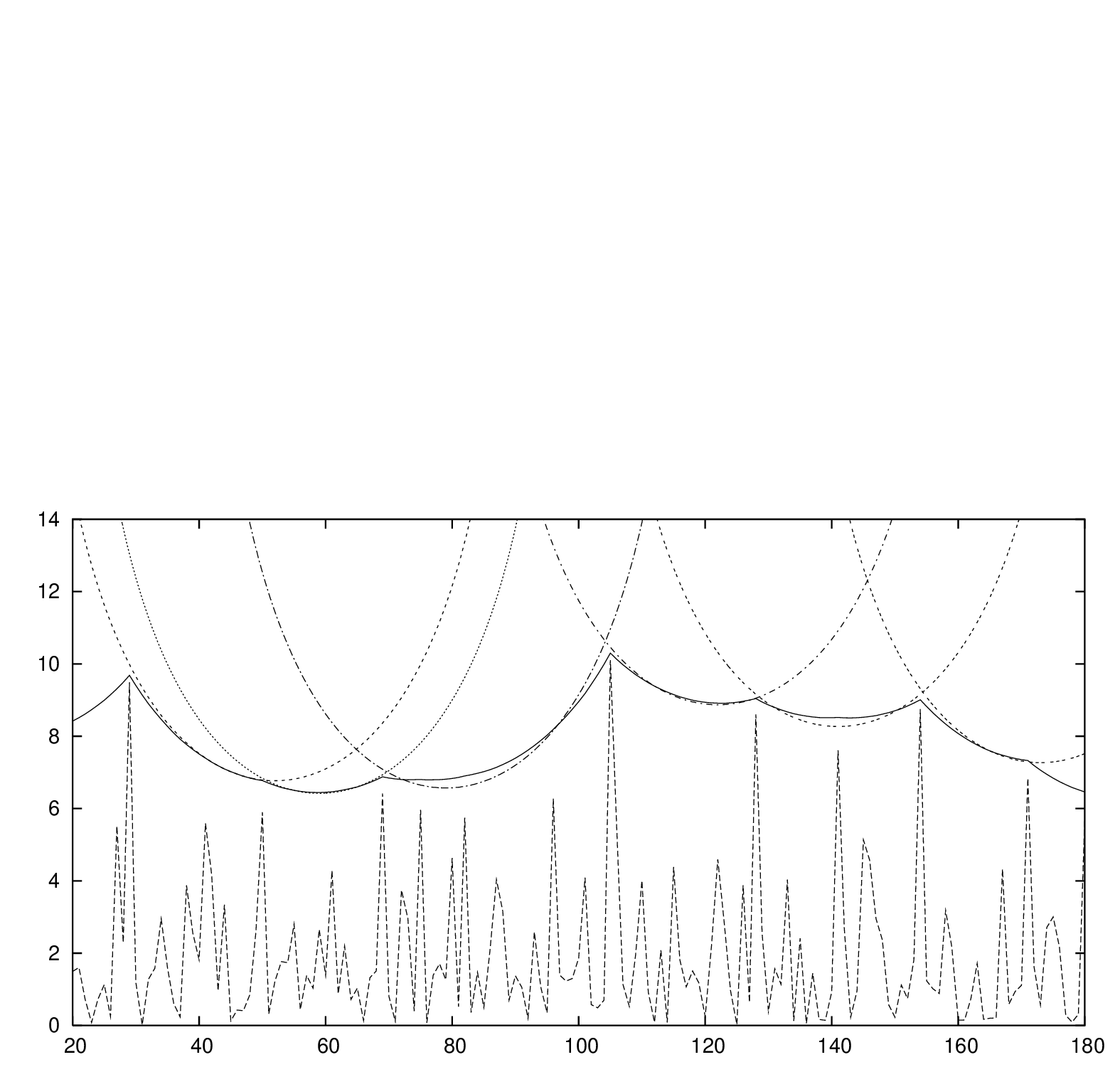}
$$
\vskip-0.5cm
\centerline{\rmsmall Fig. 2: Substrate $J_1=0$, $K_1=0.5$; film (solid line)
and Wulff shapes $J_2=5$, $K_2=0.125$.}
\medskip
The Wulff necklace is clearly a good approximation, except when the substrate
rises near but below the film surface, like near $x\simeq 80$ or $x\simeq 140$ 
on Fig. 2, which causes some entropic repulsion.
As $K_2$ decreases, the number of substrate peaks visible from the look of the
film decreases (see Fig. 1). 
The relevant peaks are large deviations events of the substrate
disorder. But there is no obvious a priori criterion to tell
which $h^1_i$ will emerge, depending upon $K_2$. The film tension leads to
correlations between the poles.

In the following sections, we simplify the model in order to understand the
selection by the film of substrate peaks as $K_2$ varies. When the substrate 
has short range correlations, its correlation length may be taken as a basic
unit of length. A natural simplification is to assume that it is also the
lattice unit, and that the substrate is i.i.d. on this scale: $J_1=0$.
Then, by scaling the unit of height, there is no loss of generality in taking
$K_1=1$, leaving two independent parameters $J_2$ and $K_2$.

The second and main simplification is to consider that when a Wulff shape
between two peaks hangs strictly above the substrate between the two peaks,
then it is not affected by the substrate between the two peaks. Some 
fluctuation effects (entropic repulsion) are neglected here, an approximation
which is better justified when $J_2$ is large.

Our aim will be to estimate the density of relevant peaks as function of $K_2$,
with different kinds of Wulff shapes. And further to get an idea of the 
probability distribution of peak localizations and heights.

\vfill\eject 
\bigskip\noindent
{\bf 3. A necklace of Wulff shapes}
\medskip\noindent
Let $h_i$, $i\in\Ze$ be an iid sequence of exponentially
distributed random variables of mean 1, and let
$$
{\S}=\{(i,z)\in\Ze\times\Re:z\le h_i\} \eqno{(3.1)}
$$
be the random set representing the substrate. The upper index 1 for the 
substrate height $h^1_i\equiv h_i$ is now omitted.

Let $1<a\le\infty$ and let $W:\ ]-a,a[\to\Re$ be a continuous even function, 
strictly increasing on $[0,a[\,$, with $W(0)=0$. If $a=\infty$, we require
$W(x)>-b+\la|x|^\al$ for some $b,\la,\al>0$, for all $x$. Examples:
$$
\eqalign{
{\rm Cone}:\qquad W(x)&=\la|x|\cr
{\rm Parabola}:\qquad W(x)&=\la x^2 \cr
{\rm Semi-circle}:\qquad W(x)&=\la^{-1}-\sqrt{\la^{-2}-x^2}\cr
{\rm SOS\ Wulff\ shape}\ \bigl(z(x)=(2.3)\bigr):\qquad 
W(x)&=z(x)-{1\over K_2}\log{2\over J_2}
} \eqno{(3.2)}
$$
with $a=\la^{-1}$ for the circle and $a=J_2/K_2$ for the SOS Wulff shape.
Translating the graph of such a function by $(x_*,h_*)\in\Re^2$ will define 
$W(x_*,h_*;\cdot)$, so that $W(0,0;\cdot)=W(\cdot)$ and
$$
W(x_*,h_*;x)=h_*+W(x-x_*) \eqno{(3.3)}
$$
The graph of this function, as a subset of $\Re^2$, is denoted $W(x_*,h_*)$.
The film over the substrate (or the coating of the substrate) is then defined 
almost surely as the graph $I$ of the function $I(x)$ defined in the following 
proposition:
\smallskip\noindent
{\bf Proposition 3.1: }{\sl Under the hypotheses on $W(\cdot)$ stated before 
(3.2), let
$$
I(x)=\inf\Bigl\{W(x_*,h_*;x):(x_*,h_*)\in\Re^2,|x_*-x|<a,\,
W(x_*,h_*)\cap\S=\emptyset\Bigr\} \eqno{(3.4)}
$$
whenever the infimum is taken over a non-empty family,
and $I(x)=+\infty$ if $W(x_*,h_*)\cap{\S}\ne\emptyset\ \forall\, (x_*,h_*)$.
Then,  almost surely, the infimum is attained and 
$I(x)<\infty\ \forall x\in\Re$.
\smallskip\noindent
Proof.} Straightforward. The hypothesis on $W(\cdot)$ could be weakened to
$W(x)>-b+\la\log|x|$ for suitable $\la$. In the following, we shall instead
strengthen the hypothesis to make $W$ convex.
\vskip-0.8cm
$$
\psboxto(13cm;0cm){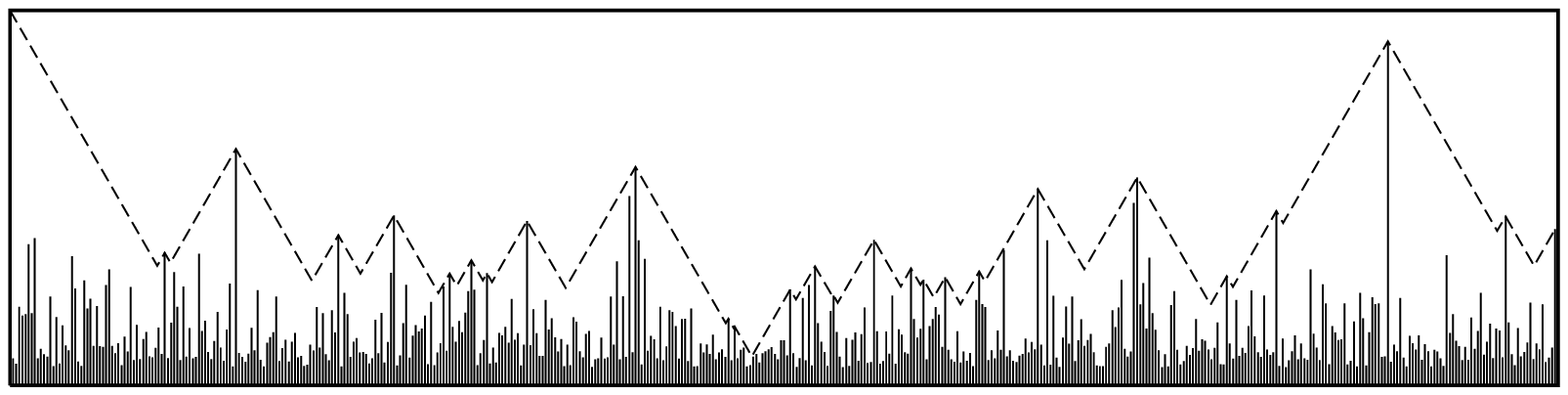}
$$
\vskip-0.5cm
\centerline{\rmsmall Fig. 3: $I(x)$ with $W(x)=\la|x|$.}
\medskip
\vskip-1cm
$$
\psboxto(13cm;0cm){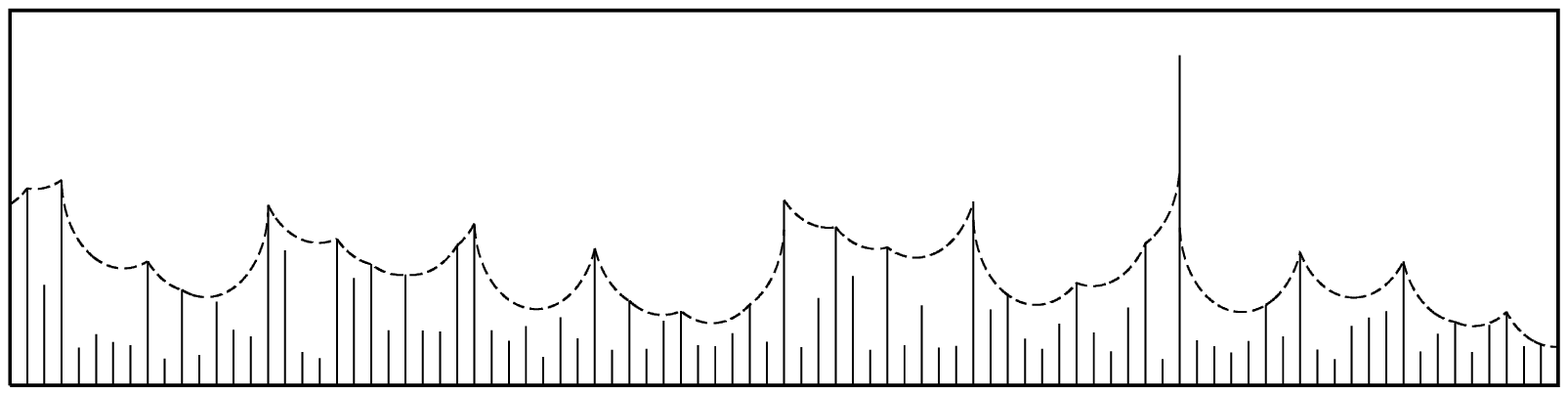}
$$
\vskip-0.5cm
\centerline{\rmsmall Fig. 4: $I(x)$ with $W(x)$ semi-circular.}
\medskip
In words: We have a ``Wulff shape'', symmetric about a vertical axis. Above 
each $x_*\in\Re$, a ``Wulff shape'' is taken down from $+\infty$, until it 
touches the substrate. The film is the infimum of the resulting collection of 
Wulff shapes. The film height $I(j)$ at an integer point $j$ models the 
thermal average of $h^2_j$ in the preceding section.

A Wulff shape in one dimension is a solution to the second 
order differential equation for the function $W(x)$,
$$
{d\over dx}{d\over d\bigl(dW/dx\bigr)}\tilde\sigma\bigl(dW/dx\bigr)=K_2 
\eqno{(3.5)}
$$
where $\tilde\sigma(\cdot)$ is the projected surface tension, or interface free
energy per unit length of interface projected onto the $x$-axis, as function of
the slope $dW/dx$. The parameter $K_2$, conjugate to the film volume 
$\sum h^2_i$ in (2.2), is the pressure, or pressure difference with the medium 
above, $\Delta p$. The surface tension  should have convexity properties such 
that the solution $W(x)$ to (3.5) is a convex function, in fact typically 
strictly convex, which we assume henceforth, except for the special 
case $W(x)=\la|x|$.

A Hamiltonian of the form (2.2) but quadratic in $(h^2_i-h^2_j)$ gives
$$
\tilde\sigma\bigl(dW/dx\bigr)=J_2\bigl(dW/dx\bigr)^2+{\rm constant}
$$
so that a solution to (3.5) in this case is $W(x)=\la x^2$ with 
$\la=K_2/(4J_2)$. Solutions to the Wulff equation (3.5) generally scale as
$$
W_{K_2}(x)=K_2^{-1}W_1(K_2x)\eqno{(3.6)}
$$
where $W_1(x)$ is a solution to (3.5) with $K_2=1$. The semi-circular shape
in (3.2) follows this scaling, with $\la$ proportional to $K_2$.

When a Wulff shape is taken down from infinity above $x=0$ until it touches the
substrate $\S$, say at $(j_0,h_{j_0})$ with $j_0\le0$ for definiteness, it 
stops as 
$$
W(0,h_*(0);x)=h_{j_0}-W(j_0)+W(x) \eqno{(3.7)}
$$
We can then let it slide to the right by an amount $x_*$, keeping contact with 
$(j_0,h_{j_0})$, as
$$
W(x_*,h_*(x_*);x)=h_{j_0}-W(j_0-x_*)+W(x-x_*) \eqno{(3.8)}
$$ 
At a given $x>0$, this is a strictly decreasing function of $x_*$ so long as 
$x_*\le x$. It is also a strictly decreasing function of $x_*$ for $x_*>x$, 
provided $W(\cdot)$ is a strictly convex function: For $y_*>x_*$,
$$
{W(y_*)-W(y_*-x)\over x}>{W(x_*)-W(x_*-x)\over x} \eqno{(3.9)}
$$
So we let the Wulff shape slide to the right keeping contact with 
$(j_0,h_{j_0})$ until it touches a second point $(k_0,h_{k_0})$:
$$
h_{k_0}-h_{j_0}=W(k_0-x_*)-W(j_0-x_*) \eqno{(3.10)}
$$ 
We thus get the unique Wulff shape going through $(j_0,h_{j_0})$ and 
$(k_0,h_{k_0})$, 
which we denote\break $W(j_0,h_{j_0},k_0,h_{k_0};x)$. If the first contact 
point is $(k_0,h_{k_0})$ with $k_0\ge0$, instead of $(j_0,h_{j_0})$ with 
$j_0\le0$, the same construction works symmetrically, sliding to the left. 
In any case we get a Wulff shape at a minimal height, so that, almost surely,
$$
I(x)=W(j_0,h_{j_0},k_0,h_{k_0};x),\qquad j_0\le x\le k_0 \eqno{(3.11)}
$$
The proposition below follows easily:
\smallskip\noindent
{\bf Proposition 3.2: }{\sl
Let $W:\Re\to\Re$  be a continuous even function, strictly convex with 
$W(0)=0$, or $W(x)=\la|x|$ with $\la>0$. Then $I(x)$ defined in (3.4)(3.1) is 
also, almost surely, 
$$
I(x)=\sup\Bigl\{W(j,h_{j},k,h_{k};x):j,k\in\Ze,\,j\le x\le k\Bigr\}
\eqno{(3.12)}
$$
Proof. } Starting from (3.4), we have $I(x)=W(x_*,h_*;x)$. This must be also
$I(x)=W(j,h_{j},k,h_{k};x)$ for some $j,k$: Otherwise (3.7-11) would give a 
smaller $I(x)$. And this $W(j,h_{j},k,h_{k};x)$ must be the supremum,
otherwise it would intersect $\S$ (as shown in more detail in the proof of
Lemma 5.1 below).
\smallskip\noindent
{\bf Remark:} For the first two examples in (3.2), we have respectively
$$
W(j,h_j,k,h_k;x)=\max\Bigl\{h_{j}-\la(x-j),\,h_{k}+\la(x-k)\Bigr\} 
\eqno{(3.13a)}
$$
$$
W(j,h_j,k,h_k;x)=\la(x-j)(x-k)+{k-x\over k-j}\,h_j
+{j-x\over j-k}\,h_k   \eqno{(3.13b)}
$$
The random interface $I(x)$ also defines, as a marginal, a point process of 
interest:
\smallskip\noindent
{\bf Proposition 3.3: }{\sl Under the same hypotheses as in Proposition 3.2, 
let
$$
\Be=I\cap{\S}=\Bigl\{i\in\Ze:I(i)=h_i\Bigr\}\eqno{(3.14)}
$$
Then
$$
\Be=\{i\in\Ze:h_i\ge W(j,h_j,k,h_k;i)\quad \forall\ j<i<k\}\eqno{(3.15)}
$$
and, almost surely, $\Be$ can be written as $\Be=\{b_n\}_{n\in\Ze}$ with 
$b_{n+1}-b_n\ge1\ \forall n$, and $b_0=\min\{b_n:b_n\ge0\}$. 
\smallskip\noindent Proof. } 
If $i$ belongs to (3.14) then Proposition 3.2 implies that it belongs 
also to (3.15). If $i$ belongs to (3.15), then we start from 
$$
h_i=W(i,h_i;i)\ge W(j,h_j,k,h_k;i)\quad \forall\ j<i<k
$$
and slide $W(i,h_i;i)$ to the right following (3.7-11) until it touches $\S$ at
$(k,h_k)$ with $h_i=W(i,h_i,k,h_k;i)$. Then Proposition 3.2 implies that $i$ 
belongs also to (3.14).

The set $\Be$ is also the set of points which can be obtained like $j_0$ or 
$k_0$ in (3.5)-(3.8), starting from any $x\in\Re$, not just $x=0$.
\bigskip\noindent
{\bf 4. Estimates}
\medskip\noindent
\smallskip\noindent
{\bf Proposition 4.1: }{\sl Let $W(x)=\la|x|$. Then
$$
\lim_{\la\searrow0}{\Pe(0\in\Be)\over\la}={1\over2}
\eqno{(4.1)}
$$
\smallskip\noindent
Proof.}
$$
\eqalign{
\Pe(0\in\Be)&=\Pe\Bigl(h_0\ge\max\bigl\{h_{j}-\la(x-j),\,h_{k}+\la(x-k)\bigr\}
\quad \forall\ j<0<k\Bigr)\cr
&=\Pe\bigl(h_0\ge h_i-\la|i|\quad\forall\,i\in\Ze_*\bigr)\cr
&=\int_0^\infty dx\,e^{-x}\prod_{i=1}^\infty\Bigl(1-e^{-x-\la i}\Bigr)^2
}\eqno{(4.2)}
$$
\noindent{---} Upper bound:
$$
\prod_{i=1}^\infty\Bigl(1-e^{-x-\la i}\Bigr)^2
<e^{-2e^{-x}\sum_{i=1}^\infty e^{-\la i}}
=e^{-2e^{-x}{e^{-\la}\over1-e^{-\la}}}\eqno{(4.3)}
$$
so that
$$
\Pe(0\in\Be)<\int_0^\infty dx\,e^{-x}e^{-2e^{-x}{e^{-\la}\over1-e^{-\la}}}
={1-e^{-\la}\over2e^{-\la}}\Bigl(1-e^{-2{e^{-\la}\over1-e^{-\la}}}\Bigr)
\eqno{(4.4)}
$$
\noindent{---} Lower bound: For $x>1$, we can use $1-\ep>e^{-\ep-\ep^2}$, with 
$\ep=e^{-x-\la i}$. Then 
$$
\prod_{i=1}^\infty\Bigl(1-e^{-x-\la i}\Bigr)^2
>e^{-2e^{-x}\sum_{i=1}^\infty e^{-\la i}-2e^{-2x}\sum_{i=1}^\infty e^{-2\la i}}
=e^{-2e^{-x}{e^{-\la}\over1-e^{-\la}}-2e^{-2x}{e^{-2\la}\over1-e^{-2\la}}}
\eqno{(4.5)}
$$
so that
$$
\eqalign{
\Pe(0\in\Be)&>\int_1^\infty dx\,e^{-x}e^{-2e^{-x}{e^{-\la}\over1-e^{-\la}}
-2e^{-2x}{e^{-2\la}\over1-e^{-2\la}}}\cr
&>\int_0^{\la\ln(\la^{-1})}dy\,e^{-2y{e^{-\la}\over1-e^{-\la}}
-2y^2{e^{-2\la}\over1-e^{-2\la}}}\cr
&>{1-e^{-\la}\over2e^{-\la}}\Bigl(1-e^{-2\la\ln(\la^{-1})
{e^{-\la}\over1-e^{-\la}}}\Bigr)
\,e^{-2(\la\ln(\la^{-1}))^2{e^{-2\la}\over1-e^{-2\la}}}
}\eqno{(4.6)}
$$
The upper and lower bounds together imply (4.1).

\smallskip\noindent
{\bf Proposition 4.2: }{\sl  Let $W(x)=\la x^2$. Then
$\exists\, a,b>0$ such that $\forall\, 0<\la<1/4$
$$
a\,{\la\over\ln(\la^{-1})}<\Pe(0\in\Be)<3\sqrt{\la\over\pi}+b\la
\eqno{(4.7)}
$$
}\smallskip\noindent
{\bf Remark:} In view of the proof of the upper bound, see below, we expect 
that $\Pe(0\in\Be)$ is of order $\sqrt\la$ for small $\la$.
\smallskip\noindent
{\sl Proof.}
$$
\Pe(0\in\Be)=\Pe\bigl(h_0\ge W(-j,h_{-j};k,h_k;0)\quad \forall\ j,k\ge1\bigr)
\eqno{(4.8)}
$$
with
$$
W(-j,h_{-j};k,h_k;0)={k\over j+k}\,h_{-j}+{j\over k+j}\,h_k-\la jk 
\eqno{(4.9)}
$$

\noindent{---} Upper bound:
$$
\eqalign{
\Pe(0\notin\Be|h_0=x)&=\Pe\bigl(\,\exists\ j,k>0:
x<W(-j,h_{-j};k,h_k;0)\,\bigr)\cr
&>\Pe\bigl(\,\exists\ j,k>0:\,h_{-j}>x+\la j^2,\,h_k>x+\la k^2\,\bigr)\cr
&=\Bigl(1-\prod_{j=1}^\infty\Pe\bigl(\,h_{-j}>x+\la j^2\bigr)\,\Bigr)^2
}\eqno{(4.10)}
$$
where we used 
$$
W\bigl(-j,x+\la j^2;k,x+\la k^2;0\bigr)>x \eqno{(4.11)}
$$ 
and $W(\cdots)$ increasing in $h_{-j}$ and in $h_k$. Then
$$
\eqalign{
\prod_{j=1}^\infty\Pe\bigl(\,h_{-j}>x+\la j^2\bigr)
=\prod_{j=1}^\infty\bigl(\,1-e^{-x-\la j^2}\,\bigr)
<e^{-e^{-x}\sum_{j=1}^\infty e^{-\la j^2}}
<e^{-e^{-x}({1\over2}\sqrt{\pi\over\la}-1)}
}\eqno{(4.12)}
$$
so that
$$
\eqalign{
\Pe(0\in\Be)
&<\int_0^\infty dx\,e^{-x}\Bigl(2e^{-e^{-x}({1\over2}\sqrt{\pi\over\la}-1)}
-e^{-2e^{-x}({1\over2}\sqrt{\pi\over\la}-1)}\Bigr)\cr
&=\int_0^1 dy\,\Bigl(2e^{-y({1\over2}\sqrt{\pi\over\la}-1)}
-e^{-2y({1\over2}\sqrt{\pi\over\la}-1)}\Bigr)\cr
&=3\sqrt{\la\over\pi}\,
{1+{1\over3}e^{2-\sqrt{\pi\over\la}}\over 1-2\sqrt{\la\over\pi}}
\qquad{\rm for}\ \la<\pi/4
}\eqno{(4.13)}
$$
\noindent{---} Lower bound: From the Harris-FKG inequality [H,FKG],
$$
\Pe\bigl(\,\forall\ j,k>0:\,x\ge W(-j,h_{-j};k,h_k;0)\,\bigr)
\ge\prod_{j=1}^\infty
\Pe\bigl(\,\forall k>0:\,x\ge W(-j,h_{-j};k,h_k;0)\,\bigr)\eqno{(4.14)}
$$
$$
\eqalign{
\Pe\bigl(\,\forall k>0:\,&x\ge W(-j,h_{-j};k,h_k;0)\,\bigr)
=\int_0^\infty dy\,e^{-y}\prod_{k=1}^\infty
\Pe\bigl(\,x\ge W(-j,y;k,h_k;0)\,\bigr)\cr
&>\int_0^{x+\la j^2}dy
\,e^{-y}\prod_{k=1}^\infty\Bigl(1-e^{-\la k^2-{k\over j}(x+\la j^2-y)-x}\Bigr)
\cr
>\int_0^{x+\la j^2}dy\,
&\exp\Bigl(-y-\sum_{k=1}^\infty e^{-\la k^2-{k\over j}(x+\la j^2-y)-x}
-\sum_{k=1}^\infty e^{-2\la k^2-{2k\over j}(x+\la j^2-y)-2x}\Bigr)
}\eqno{(4.15)}
$$
where we used $1-\ep>e^{-\ep-\ep^2}$, with $\ep<e^{-x}$ and $x\ge1$ henceforth.
Then
$$
\sum_{k=1}^\infty e^{-\la k^2-{k\over j}(x+\la j^2-y)-x}
<\sum_{k=1}^\infty e^{-{k\over j}(x+\la j^2-y)-x}
={e^{-x}\over e^{x+\la j^2-y\over j}-1}
<{je^{-x}\over x+\la j^2-y}\eqno{(4.16)}
$$
$$
\eqalign{
\Pe\bigl(\,\forall k>0:\,x\ge W(-j,h_{-j};k,h_k;0)\,\bigr)
&>\int_0^{x+\la j^2\over2}dy\,e^{-y-{je^{-x}\over x+\la j^2-y}}\cr
&>{e^{-{je^{-x}\over x+\la j^2}}\over1+{2je^{-x}\over(x+\la j^2)^2}}
\Bigl(1-e^{-{x+\la j^2\over2}-{je^{-x}\over x+\la j^2}}\Bigr)\cr
>\exp\Bigl(-{je^{-x}\over x+\la j^2}-{2je^{-x}\over(x+\la j^2)^2}
&-e^{-{x+\la j^2\over2}-{je^{-x}\over x+\la j^2}}
-e^{-(x+\la j^2)-{2je^{-x}\over x+\la j^2}}\Bigr)
}\eqno{(4.17)}
$$
where $(1-\al y)^{-1}<1+2\al y$ for $0<\al y<1/2$ has been used before 
integration over $y$, and then again $1-\ep>e^{-\ep-\ep^2}$ and also
$1/(1+X)>e^{-X}$. Now
$$
\eqalign{
\prod_{j=1}^{\la^{-1}}
\Pe\bigl(\,\forall k>0:\,x\ge W(-j,&h_{-j};k,h_k;0)\,\bigr)
>\exp\Bigl(-{e^{-x}\over2\la}\ln{x+\la^{-1}\over x+\la}
-{e^{-x}\over2\sqrt{\la x}}\cr
&-{e^{-x}\over2\la(x+\la)}-{e^{-x}\over\sqrt{\la x^3}}
-\sqrt{2\pi\over\la}e^{-x/2}-\sqrt{\pi\over\la}e^{-x}-2\Bigr)
}\eqno{(4.18)}
$$
The sums over $j\in\Ze_+$ in the exponential were bounded by corresponding 
integrals over $[1,\la^{-1}]$ or $\Re_+$ plus a bound of the maximum of the 
integrand. Then with $\la<1/4$ and $x>4$,
$$
\eqalign{
\prod_{j=1}^{\la^{-1}}
\Pe\bigl(\,\forall k>0:\,x\ge W(-j,h_{-j};k,h_k;0)\,\bigr)
>\exp\Bigl(-{e^{-x}\over2\la}\ln{\la^{-1}}
-{e^{-x}\over\la}-{\sqrt{2\pi}e^{-{x\over2}}\over\sqrt\la}-2\Bigr)
}\eqno{(4.19)}
$$
For $j>\la^{-1}$, the range of integration in (4.15) is chosen as 
$0<y<x+\la j^2-j\ln 2$. Then only the first few $k=1,2,\dots$ play a role for 
the event $\exists\,k:x< W(-j,y;k,h_k;0)$. Precisely:
$$
\sum_{k=1}^\infty e^{-\la k^2-{k\over j}(x+\la j^2-y)-x}
<{e^{-x}\over e^{x+\la j^2-y\over j}-1}
<e^{-x-{x+\la j^2-y\over j}}+2e^{-x-2{x+\la j^2-y\over j}}\eqno{(4.20)}
$$
$$
\eqalign{
\int_0^{x+\la j^2-j\ln 2}dy\,&\exp\Bigl(-y
-e^{-x-{x+\la j^2-y\over j}}-2e^{-x-2{x+\la j^2-y\over j}}\Bigr)\cr
&>\int_0^{x+\la j^2-j\ln 2}dy\,e^{-y}
\Bigl(1-3e^{-x-{x+\la j^2-y\over j}}\Bigr)\cr
&=1-e^{-(x+\la j^2-j\ln 2)}-3{e^{-x-{x+\la j^2\over j}}\over1-1/j}
\Bigl(1-e^{-(x+\la j^2-j\ln 2)(1-1/j)}\Bigr)\cr
&>1-5e^{-x-\la j}>\exp(-5\,e^{-x-\la j}-25\,e^{-2x-2\la j})\cr
}\eqno{(4.21)}
$$
for $\la<1/4$ and $x>2$. Then
$$
\eqalign{
\prod_{j>\la^{-1}}
\Pe\bigl(\,\forall k>0:\,x\ge W(-j,h_{-j};k,h_k;0)\,\bigr)
&>\prod_{j>\la^{-1}}\exp(-5\,e^{-x-\la j}-25\,e^{-2x-2\la j})\cr
&\ge\exp\Bigl(-{5e^{-x-1}\over1-e^{-\la}}-{25e^{-2x-2}\over1-e^{-2\la}}\Bigr)
}\eqno{(4.22)}
$$
Putting together (4.19) and (4.22) and integrating with $d(e^{-x})$ over an 
interval
$$
\al\,{\la\over\ln(\la^{-1})}<e^{-x}<\beta\,{\la\over\ln(\la^{-1})}\eqno{(4.23)}
$$
yields the lower bound and completes the proof of the proposition.
\bigskip\noindent
{\bf 5. Gibbs measure}
\medskip\noindent
So far the equivalent definitions of $\Be$, through (3.14) or (3.15), require
a knowledge of the whole system in order to decide whether a point $i\in\Be$.
Here we will define local Gibbs factors where $(b_n,h_{b_n})$ is coupled to
$(b_{n-1},h_{b_{n-1}})$  and $(b_{n+1},h_{b_{n+1}})$ only. We first derive an
equivalent definition of $\Be$:

\smallskip\noindent
{\bf Proposition 5.1: }{\sl 
$\Be$ defined in Proposition (3.3) is almost surely the only countable ordered 
subset of $\Ze$, denoted $\Be=\{b_n\}_{n\in\Ze}$, 
with $b_0=\min\{b_n\ge0\}$, obeying the following two conditions:
$$
h_i<W(b_n,h_{b_n},b_{n+1},h_{b_{n+1}};i)\qquad\forall\quad 
b_n< i< b_{n+1} \eqno{(5.1)}
$$
$$
h_{b_n}\ge W(b_{n-1},h_{b_{n-1}},b_{n+1},h_{b_{n+1}};b_n)\qquad\forall\quad
 n\in\Ze \eqno{(5.2)}
$$
$\Be$ is  the minimal subset of $\Ze$ such that the collection of Wulff
shapes suspended from this subset lies above the rest of the substrate, i.e.
the minimal subset of $\Ze$ satisfying (5.1). 
\smallskip\noindent Proof. }
From Propositions (3.2) and (3.3), $\Be$ is almost surely the only countable 
ordered subset of $\Ze$ satisfying (5.1) and  
$$
h_{b_n}\ge W(b_{n-p},h_{b_{n-p}},b_{n+q},h_{b_{n+q}};b_n)\,,\qquad\forall\quad
 n\in\Ze,\,p\ge1,q\ge1 \eqno{(5.2')}
$$
We only need to prove that (5.2) implies (5.2'). 
Let us assume two instances of (5.2),
$$
h_{b_1}\ge W(h_{b_0},h_{b_2};b_1)\eqno{(5.3a)}
$$
$$
h_{b_2}\ge W(h_{b_1},h_{b_3};b_2)\eqno{(5.3b)}
$$
and, contradicting (5.2') with $p=2,q=1$,
$$
h_{b_2}<W(h_{b_0},h_{b_3};b_2)\qquad ??? \eqno{(5.3')}
$$
Let us show that $(5.3a-b)$ with the absurd $(5.3')$ would imply that 
$W(h_{b_1},h_{b_3};\cdot)$ and $W(h_{b_0},h_{b_2};\cdot)$ have two
intersections. Indeed 
$$
\eqalign{
W(h_{b_1},h_{b_3};b_1)&\ge W(h_{b_0},h_{b_2};b_1)\qquad{\rm using}\
 (5.3a)\cr
W(h_{b_1},h_{b_3};b_2)&\le W(h_{b_0},h_{b_2};b_2)\qquad{\rm using}\
 (5.3b)\cr
W(h_{b_1},h_{b_3};b_3)&=W(h_{b_0},h_{b_3};b_3)\cr
&=W(h_{b_0},W(h_{b_0},h_{b_3};b_2);b_3)\cr
&\ge W(h_{b_0},h_{b_2};b_3)\qquad{\rm using}\ (5.3')
}
$$
so that $W(h_{b_1},h_{b_3};\cdot)$ is higher than 
$W(h_{b_0},h_{b_2};\cdot)$ at $b_1$ and $b_3$ and lower at $b_2$, which
implies two intersections, impossible for two Wulff shapes which are translates
of one another.
Therefore $(5.3')$ cannot hold true; (5.2) with $p=1,q=1$ implies (5.2) also 
with $p=2,q=1$. The argument extends easily to all $p,q$ and the proof of
Proposition (5.1) is readily completed.

Let $l_n=b_n-b_{n-1}$ and $x_n=h_{b_n}$. A Gibbs measure formulation
for $\{l_n,x_n\}_{n\in\Ze}$ starts from i.i.d. a priori measures: 
Counting measure on $\Ze_+$ for each $l_n$, exponential distribution 
$\exp(-x_n)dx_n$ on $\Re_+$ for each $x_n$. And a product of
local Gibbs factors,
$$
\prod_n F(x_n,l_{n+1},x_{n+1})
\prod_n G(x_{n-1},l_n,x_n,l_{n+1},x_{n+1}) \eqno{(5.4)}
$$
with
$$
F(x_0,l_1,x_1)=\prod_{i=1}^{l_1-1}\Bigl(1-e^{-W(0,x_0,l_1,x_1;i)}\Bigr)
\eqno{(5.5)}
$$
and
$$
G(x_{-1},l_0,x_0,l_1,x_1)=
1_{x_0\,\ge\, W(-l_0,x_{-1},l_1,x_1;0)} \eqno{(5.6)}
$$

Let us define explicitly a finite volume Gibbs measure with free boundary 
conditions: 
\smallskip\noindent
{\bf Proposition 5.2: }{\sl Let $h_i$, $i\in\Ze$ be an iid sequence of 
exponentially distributed random variables of mean 1. Let $W:\Re\to\Re$  be a 
continuous even function, strictly convex with $W(0)=0$. Let
$$
I_{[0,L]}(i)=\sup\Bigl\{W(j,h_{j},k,h_{k};i):\,0\le j\le i\le k\le L\Bigr\} 
\eqno{(5.7)}
$$
and
$$
\Be_{[0,L]}=\Bigl\{i\in[0,L]\cap\Ze:I_{[0,L]}(i)=h_i\Bigr\}
=\{b_0,\dots,b_N\}\eqno{(5.8)}
$$
with $b_0=0$, $b_N=L$, and $N\ge1$. Let  $l_n=b_n-b_{n-1}$ and $x_n=h_{b_n}$. 
Then \break
$\Bigl\{N,\{l_1,\dots,l_N\},\{x_0,\dots,x_N\}\Bigr\}$ is distributed according 
to
$$
\Xi_{[0,L]}^{-1}\prod_0^N dx_ne^{-x_n}
\prod_0^{N-1} F(x_n,l_{n+1},x_{n+1})
\prod_1^{N-1} G(x_{n-1},l_n,x_n,l_{n+1},x_{n+1}) \eqno{(5.9)}
$$
where the partition function is
$$
\Xi_{[0,L]}=\sum_{N=1}^{L}\sum_{l_1\dots l_N}
\int\prod_0^N dx_ne^{-x_n}
\prod_0^{N-1} F(x_n,l_{n+1},x_{n+1})
\prod_1^{N-1} G(x_{n-1},l_n,x_n,l_{n+1},x_{n+1}) \eqno{(5.10)}
$$
and the sum over the positive integers $l_1\dots l_N$ is constrained by
$l_1+\dots+ l_N=L$.
\smallskip\noindent
Proof.}
Proposition 5.2 is a simple corollary of Proposition 5.1. 

Such a model is solvable in principle. A natural first step is to change to
a pressure ensemble with $\exp(-pL)$ and $L$ random in order to get rid of
the global constraint over $l_1\dots l_N$.

\bigskip\noindent
{\bf Acknowledgments}: F. Dunlop acknowledges efficient support and kind
hospitality at Universit\'e de Mons-Hainaut and CRMM where much of the present 
work was done.
\medskip\noindent

\bigskip\noindent
{\bf References}
\medskip\noindent

\noindent\item{[AS1]} D.B. Abraham, E.R. Smith: {\sl Surface film thickening: 
An exactly solved model}, Phys. Rev. {\bf B26}, 1480--1482, (1982).

\noindent\item{[AS2]} D.B. Abraham, E.R. Smith: {\sl An exactly solved model
with a wetting transition}, J. Stat. Phys. {\bf 43}, 621--643 (1986).

\noindent\item{[BN]} H. van Beijeren, I. Nolden: pp 259--300 in Structure and 
Dynamics of Surfaces II, edited by W. Schommers and P. von Blanckenhagen.
Topics in Current Physics Vol. 43 (Springer-Verlag, Berlin Heidelberg, 1987).

\noindent\item{[DD]} J. De Coninck, F. Dunlop: {\sl Partial to complete 
wetting: A microscopic derivation of the Young relation}, J. Stat. Phys. 
{\bf 47}, 827--849 (1987).

\noindent\item{[FKG]} C. Fortuin, P. Kasteleyn, J. Ginibre: {\sl Correlation
inequalities on some partially ordered sets}, Comm. Math. Phys. {\bf 79},
141--151 (1981).

\noindent\item{[H]} T. E. Harris: {\sl A lower bound for the critical 
probability in a certain percolation process}, Proc. Camb. Phil. Soc. {\bf 59},
13--20 (1960).

\bye